\documentclass[preprint]{elsarticle}
\usepackage{amsmath,amssymb,a4wide}
\usepackage[curve,matrix,arrow,cmtip]{xy}
\usepackage{hyperref}

\journal{Journal of Number Theory}

\def\Ksep{{K^{\rm sep}}}

 \def\Aut{\mathop{\rm Aut}\nolimits}

 \def\Gal{\mathop{\rm Gal}\nolimits}
 \def\End{\mathop{\rm End}\nolimits}

 \def\Spec{\mathop{\rm Spec}\nolimits}
 
 \def\deg{\mathop{\rm deg}\nolimits}

\def\Quot{\mathop{\rm Quot}\nolimits}

\def\diag{\mathop{\rm diag}\nolimits}

\def\rank{\mathop{\rm rank}\nolimits}

\def\GL{\mathop{\rm GL}\nolimits}
\def\SL{\mathop{\rm SL}\nolimits}

\def\PGL{\mathop{\rm PGL}\nolimits}
\def\PSL{\mathop{\rm PSL}\nolimits}

\def\der{{\rm der}}

\def\sep{{\rm sep}}

\let\phi\varphi
\let\epsilon\varepsilon

\newtheorem{Thm}{Theorem}
\newtheorem{Prop}[Thm]{Proposition}
\newtheorem{Lem}[Thm]{Lemma}

\newtheorem{Cor}[Thm]{Corollary}

\newtheorem{Def}[Thm]{Definition}

\newtheorem{Rem}[Thm]{Remark}

\def\qed{{\hskip0pt\unskip\unskip\nobreak\hfil\penalty50
          \hskip1em\hbox{}\nobreak\hfil
           {$\square$}
          \parfillskip=0pt\finalhyphendemerits=0
          \par}\medskip}

\newenvironment{Proof}
               {\noindent{\bf Proof.}\ }
               {\qed}

\newenvironment{Proofof}[1]
               {\noindent{\bf Proof of #1.}\ }
               {\qed}

\newcommand{\BC}{{\mathbb{C}}}

\newcommand{\BF}{{\mathbb{F}}}
\newcommand{\BG}{{\mathbb{G}}}

\newcommand{\BN}{{\mathbb{N}}}

\newcommand{\BP}{{\mathbb{P}}}
\newcommand{\BQ}{{\mathbb{Q}}}

\newcommand{\BZ}{{\mathbb{Z}}}

\newcommand{\CF}{{\cal F}}

\newbox\mybox
\def\arrover#1{\mathrel{
       \setbox\mybox=\hbox spread 1.4em
              {\hfil$\scriptstyle#1$\hfil}
       \vbox{\offinterlineskip\copy\mybox
             \hbox to\wd\mybox{\rightarrowfill}}}}

\def\larrover#1{\mathrel{
       \setbox\mybox=\hbox spread 1.4em
              {\hfil$\scriptstyle#1\vphantom{g}$\hfil}
       \vbox{\offinterlineskip\copy\mybox
             \hbox to\wd\mybox{\leftarrowfill}}}}

\def\ontoover#1{\mathrel{
       \setbox\mybox=\hbox spread 1.4em
              {\hfil$\scriptstyle#1\vphantom{g}$\hfil}
       \vbox{\offinterlineskip\copy\mybox
             \hbox to\wd\mybox{\rightarrowfill\hskip-2.8mm
                               $\rightarrow$}}}}
\def\leftontoover#1{\mathrel{
       \setbox\mybox=\hbox spread 1.4em
              {\hfil$\scriptstyle#1\vphantom{g}$\hfil}
       \vbox{\offinterlineskip\copy\mybox
             \hbox to\wd\mybox{$\leftarrow$\hskip-2.8mm
                               \leftarrowfill}}}}
\let\longto\longrightarrow
\let\into\hookrightarrow

\def\Cinf{{\BC}_\infty}

\def\kinf{k_\infty}

\renewcommand{\Im}{{\mathrm {Im}}}
\def\Vo{V'}
\def\Wo{W'}

\begin{document}

\begin{frontmatter}

\title{Explicit Drinfeld Moduli Schemes and \\ Abhyankar's Generalized Iteration Conjecture}
\author{Florian Breuer}
\address{Stellenbosch University, Stellenbosch, South Africa\\ 
{\em fbreuer@sun.ac.za}\\
 {Dedicated to Mira Breuer on the occasion of her 0th birthday}}

\begin{abstract}
Let $k$ be a field containing $\BF_q$. Let $\psi$ be a rank $r$ Drinfeld $\BF_q[t]$-module determined by 
$\psi_t(X) = tX+a_1X^q+\cdots+a_{r-1}X^{q^{r-1}}+X^{q^r}$, where $t,a_1,\ldots,a_{r-1}$ are algebraically independent over $k$. Let $n\in\BF_q[T]$ be a monic polynomial. We show that the Galois group of $\psi_n(X)$ over $k(t,a_1,\ldots,a_{r-1})$ is isomorphic to $\GL_r(\BF_q[t]/n\BF_q[t])$, settling a conjecture of Abhyankar. Along the way we obtain an explicit construction of Drinfeld moduli schemes of level $tn$.
\end{abstract}


\end{frontmatter}
\section{Introduction}

A classical theorem of Weber (see \cite[Chapter 6, Corollary 1]{LangEF}) states that, if $E$ is an elliptic curve over $K=\BQ(j)$ with transcendental $j$-invariant $j$, and $K_n = K(E[n])$ denotes the field obtained by adjoining the coordinates of all $n$-torsion points of $E$ to $K$, then 
\[
\Gal(K_n/K) \cong \GL_2(\BZ/n\BZ).
\]
The goal of this paper is to prove the analogous statement for Drinfeld modules. Our result was conjectured by S.~S.~Abhyankar in \cite[\S19]{Ab01}, who didn't know about Drinfeld modules at the time, and called it the ``Generalized Iteration Conjecture''. It was part of his quest to find nice equations for nice groups.

Our approach is based on the observation that a particular morphism of Drinfeld moduli schemes is \'etale with a suitable Galois group, a result due to V.~G.~Drinfel'd \cite{Drinfeld}. 

We start by constructing some suitable rings in \S\ref{Sec:rings}, and then in \S\ref{Sec:moduli} we prove our first main result, Theorem \ref{moduli}, which gives an explicit construction of the Drinfeld moduli scheme of level $tn$. This is a generalization of the level $t$ construction due to R.~Pink \cite{Pink}.

In \S\ref{Sec:generic} we prove our second main result, Theorem \ref{generic}, which shows that this moduli scheme can be obtained from the torsion module of any ``sufficiently generic'' Drinfeld module, including the one defined by Abhyankar. Both of these results may be of independent interest.

In \S\ref{Sec:Abhyankar} we then state the third main result, Theorem \ref{Abhyankar}, which settles Abyhankar's Generalized Iteration Conjecture. This is proved in \S\S\ref{Sec:MooreCarlitz}--\ref{Sec:proof}.

Most of this paper will be comprehensible to anybody familiar with the basics of Drinfeld modules over a field, see for example \cite[Chapter 4]{GossBS} or~\cite[Chapter 12]{Rosen}. To state Theorem~\ref{moduli}, we need Drinfeld modules over a scheme, for which we recommend the exposition in \cite{Lehmkuhl}.

\paragraph{Acknowledgements} The author would like to thank Dirk Basson for helpful discussions, and the anonymous referee for useful comments and corrections. The author is particularly grateful to another anonymous referee who read an earlier draft of this paper, in which a more complicated proof of Theorem \ref{Abhyankar} was proposed, and whose suggestions lead to the approach presented here. This work was supported by NRF grant number  BS2008100900027.

%
%
%



\section{Rings generated by torsion points}\label{Sec:rings}

Let $\BF_q$ denote a finite field of $q$ elements, where $q$ is a power of the prime $p$, and let $V\cong \BF_q^r$ denote an $\BF_q$-vector space of dimension $r\geq 1$. We denote by $\Vo=V\smallsetminus\{0\}$ the set of non-zero vectors in $V$.

Denote by $A=\BF_q[t]$ the polynomial ring over $\BF_q$, let $n\in A$ be a monic polynomial and set $B:=\BF_q[t,\frac{1}{tn}] = A[\frac{1}{tn}]$.

We denote by $S_V=B[v \;|\; v\in\Vo]$ the symmetric algebra of $V$ over $B$, which is isomorphic to a polynomial ring over $B$ in $r$ independent variables; by $K_V$ the quotient field of $S_V$; by $R_V=B[\frac{1}{v}\;|\;v\in\Vo]$ the $B$-subalgebra of $K_V$ generated by $\frac{1}{v}$ for every $v\in\Vo$; and finally by $RS_V=B[v,\frac{1}{v}\;|\;v\in\Vo]$ the $B$-subalgebra of $K_V$ generated by $S_V$ and $R_V$. 
The rings $S_V$, $R_V$ and $RS_V$ are $\BZ$-graded $B$-algebras with respect to the grading $\deg(v)=-1$ and $\deg(\frac{1}{v})=1$ for all $v\in\Vo$. 
The homogeneous component of $RS_V$ of degree zero is denoted $RS_{V,0}$. 
We have
\[
RS_{V,0} = B\left[\left.\frac{v}{v'} \; \right|\; v,v'\in\Vo\right].
\]
%
These definitions essentially come from \cite{PinkSchieder}.

We define a rank $r$ Drinfeld $A$-module $\phi$ over $K_V$ by setting
\begin{eqnarray*}
\phi_t(X) & := & tX\prod_{v\in\Vo}\left(1-\frac{1}{v}X\right) \\
& = & tX + g_1X^q + \cdots + g_rX^{q^r} \in R_V[X].
\end{eqnarray*}

The coefficients $g_1,\ldots,g_r$ are algebraically independent over $B$, since $K_V$ has transcendence degree $r+1$ over $\BF_q$ while being finite over $\BF_q(t,g_1,\ldots,g_r)$.
The highest coefficent $g_r=t\prod_{v\in\Vo}\frac{1}{v}$ is a unit in $RS_V$, and by construction the $t$-torsion submodule of $\phi$ is $\phi[t]=V\subset K_V$.

Next, we want to define similar rings for the $tn$-torsion module of $\phi$.
We construct the ring $RS_W$ via generators and relations, as follows. Choose a basis $v_1,v_2,\ldots,v_r$ of $V$. 

Inside the polynomial ring $RS_V[w_1,w_2,\ldots,w_r]$, in $r$ independent variables, we consider the ideal 
%
\[
I_W := \langle \phi_n(w_1)-v_1, \; \phi_n(w_2)-v_2, \ldots, \phi_n(w_r)-v_r \rangle. 
\]
Then we define 
\[
RS_W := RS_V[w_1,w_2,\ldots,w_r]/I_W.
\] 

Consider the following subset of $RS_W$:
\[
W:=\{\phi_{a_1}(w_1) + \phi_{a_2}(w_2) + \cdots + \phi_{a_r}(w_r) \;|\; a_1,\ldots,a_r\in A/tnA\} \subset RS_{W}.
\]
By abuse of notation, we have written $\phi_a(w_i)$ with $a\in A/tnA$ when we actually mean that $a\in A$ represents a certain class in $A/tnA$. Since $\phi_a(w_i)\equiv \phi_b(w_i) \pmod{\phi_n(w_i)-v_i}$ if $a\equiv b \pmod{tn}$, it does not matter which representative we take.

\begin{Prop}\label{Prop:rings} $\;$
\begin{enumerate}
\item[(i)] The inclusion $RS_V\into RS_W$ is \'etale and $RS_W$ is reduced.
\item[(ii)] The action $a\cdot w := \phi_a(w)$, for $a\in A$ and $w\in W$, turns $W$ into an $A/tnA$-module, which is free of rank $r$, and $V\subset W$.
\item[(iii)] Every element of $\Wo:=W\smallsetminus \{0\}$ is invertible in $RS_W$.
\item[(iv)] The following identity holds in $RS_W[X]$:
\[
\phi_{tn}(X) = tnX\prod_{w\in\Wo}\left(1-\frac{1}{w}X\right).
\]
\end{enumerate}
\end{Prop}

\begin{Proof}
Since $\displaystyle\det\left(\frac{\partial}{\partial w_i}(\phi_n(w_j)-v_j)\right) = n^r \in RS_W^*$, it follows from \cite[Tag 03PC]{StacksProject} that $RS_V \into RS_W$ is \'etale. Since $RS_V$ is reduced, so is $RS_W$, by [{\em loc. cit.}], proving (i).

Every prime ideal of $RS_W$ is the kernel of a ring homomorphism $\theta : RS_W \longto F$, where $F$ is an algebraically closed field. To such a $\theta$ we associate the Drinfeld module $\phi^\theta$ over $F$ by 
\[
\phi^\theta_t(X) := \theta(t)X\prod_{v\in\Vo}\left(1-\theta\left(\frac{1}{v}\right)X\right)\in F[X].
\]
Each $w_i$ maps to $\theta(w_i)\in L$ satisfying $\phi^\theta_{n}(\theta(w_i)) = \theta(v_i)$. Since the characteristic $\ker(\theta|_A)$ of $\phi^\theta$ is prime to $tn$, the $\theta(v_1),\ldots,\theta(v_r)$ generate $\phi^\theta[t] \cong (A/tA)^r$, and thus also the $\theta(w_1),\ldots, \theta(w_r)$ generate $\phi^\theta[tn] \cong (A/tnA)^r$. From this follows that $\theta$ maps the set $W$ isomorphically to $\phi^\theta[tn]$, completing the proof of (ii).

Furthermore, $W \cap \ker(\theta) = \{0\}$ for every such $\theta : RS_W \to F$, so $W'\subset RS_W^*$, proving (iii).

%
%
%
Lastly, for each $\theta : RS_W \to F$ we have 
\[
\phi^\theta_{tn}(X) = \theta(tn)X\prod_{w\in\Wo}\left(1-\theta\left(\frac{1}{w}\right)X\right)\in F[X],
\]
since both polynomials have the same roots and linear term. It follows that each coefficient of 
\[
\phi_{tn}(X) - tnX\prod_{w\in\Wo}\left(1-\frac{1}{w}X\right)\in RS_W[X]
\]
lies in $\cap_\theta\ker(\theta) = \{0\}$, since $RS_W$ is reduced, which completes the proof of (iv).
%
%
%
%
%
%
\end{Proof}

We note that $RS_W$ is generated over $RS_V$ by the elements of $W$. 
At this point it is far from clear that $RS_W$ is integral, but this will be shown later (Theorem~\ref{Drinfeld}).

Next, we define a ring generated only by the quotients of torsion points.

Recall that $v_1\in\Vo$ is fixed. The Drinfeld module $\phi':=v_1^{-1}\phi v_1$, defined by 
\[
\phi'_t(X) = v_1^{-1}\phi_t(v_1X) = tX\prod_{v\in\Vo}\left(1-\frac{v_1}{v}X\right) \in RS_{V,0}[X],
\]
is isomorphic to $\phi$ over $K_V$.

Inside the polynomial ring $RS_{V,0}[w'_1,w'_2,\ldots,w'_r]$, in $r$ independent variables, we define the ideal 
\[
I_{W,0} := \left\langle \phi'_n(w'_1)-\frac{v_1}{v_1}, \; \phi'_n(w'_2) - \frac{v_2}{v_1}, \ldots,  \phi'_n(w'_r) - \frac{v_r}{v_1} \right\rangle. 
\]
Then we define
\[
RS_{W,0} := RS_{V,0}[w'_1,w'_2,\ldots,w'_r]/I_{W,0}.
\]
The ring $RS_{W,0}$ embeds into $RS_W$ via $w'_i \mapsto w_i/v_1$, and the above relations reflect the fact that $\phi'_n(w_i/v_1) = v_i/v_1$ for $i=1,2,\ldots,r$.
Thus we have
\[
RS_{W,0} = RS_{V,0}\left[\left.\frac{w}{w'} \;\right|\; w,w'\in\Wo\right] = B\left[\left.\frac{w}{w'} \;\right|\; w,w'\in\Wo\right] \subset RS_W.
\]
Moreover, we see that $RS_{W,0}$ is the degree zero component of $RS_W$ with respect to the $\BZ$-grading defined by $\deg(w)=-1$ for all $w\in\Wo$. 

Lastly, it follows from Proposition \ref{Prop:rings}.(iv) that 
\[
\phi'_{tn}(X) = tnX\prod_{w\in\Wo}\left(1-\frac{v_1}{w}X\right) \in RS_{W,0}[X].
\]

%


\section{Explicit Drinfeld moduli schemes}\label{Sec:moduli}

Let $S$ be a scheme over $\Spec B$. Then recall (see e.g. \cite{Lehmkuhl}) that a Drinfeld $A$-module of rank $r$ over $S$ is a pair $(L,\psi)$, where $L$ is the additive group scheme of a line bundle over $S$, and 
\[
\psi : A \longto \End_{\BF_q}(L), \qquad a\longmapsto \psi_a,
\] 
is a ring homomorphism that is defined over a trivializing open $\Spec(R)\subset S$ by 
\[
t \longmapsto \psi_t(X) = tX + e_1X^q + 
\cdots + e_nX^{q^n},
\]
where for each $i=1,2,\ldots,n$ we have $e_i\in R$, $e_r\in R^*$ and $e_i$ is nilpotent for all $i>r$. We usually drop the $L$ from our notation and refer to the Drinfeld module as $\psi$.
When $e_i=0$ for all $i>r$, we say that the Drinfeld module $\psi$ is {\em standard}. 

The $tn$-torsion submodule $\psi[tn]$ of $\psi$ is the closed subscheme of $L$ defined locally over $\Spec R$ by $\Spec \big(R[X]/\langle \psi_{tn}(X)\rangle\big)$. When $R=F$ is a field, we identify $\psi[tn]$ with $\{x\in\bar{F} \;|\; \psi_{tn}(x)=0\}$.

A level-$tn$ structure on $\psi$ is a homomorphism of $A$-modules
\[
\mu : \big((tn)^{-1}A/A\big)^r \longto L(S)
\]
which induces an equality of divisors 
\[
\sum_{\alpha\in \big((tn)^{-1}A/A\big)^r} \mu(\alpha) = \psi[tn].
\]

For Drinfeld modules over $\Spec R$, where $R$ is a $B$-algebra, this is equivalent to the following more convenient formulation. 
Fix $A$-module isomorphisms 
\[
((tn)^{-1}A/A)^r\stackrel{\sim}{\longto} (A/tnA)^r\stackrel{\sim}{\longto} W\quad\text{and}\quad (t^{-1}A/A)^r\stackrel{\sim}{\longto}(A/tA)^r\stackrel{\sim}{\longto} V
\] 
such that the following diagram commutes: 
\[
\xymatrix{
\big((tn)^{-1}A/A\big)^r\ar[r]^{\sim} & \big(A/tnA\big)^r\ar[r]^{\sim} & W \\
\big(t^{-1}A/A\big)^r\ar[r]^{\sim}\ar@{^{(}->}[u] & \big(A/tA\big)^r\ar[r]^{\sim} & V\ar@{^{(}->}[u] \\
}
\]

Then a level-$tn$ structure on a Drinfeld module $\psi$ over $\Spec R$ is equivalent to an $A$-module homomorphism (where the $A$-module structure on $R$ is induced by $\psi$)
\[
\mu : W \longto R
\]
such that $\mu(\Wo)\subset R^*$ and
\[
\psi_{tn}(X) = tnX\prod_{w\in\Wo}\left(1-\frac{X}{\mu(w)}\right) \in R[X].
\]

Here we have made essential use of the fact that the characteristic $\ker(A \to R)$ of $\psi$ is prime to $tn$, since $R$ is a $B$-algebra and $tn\in B^*$.
%

In particular, by Proposition \ref{Prop:rings}, $\phi'$ carries the level-$tn$ structure
\[
\lambda : W \longto RS_{W,0}; \quad w\mapsto \frac{w}{v_1}.
\]

Our first main result is the fact that $\Spec(RS_{W,0})$ is the fine moduli scheme for rank $r$ Drinfeld $A$-modules with level-$tn$ structure. 
Denote by $E=\BG_{\mathrm{a},RS_{W,0}}$ the additive group scheme over $\Spec(RS_{W,0})$. Then the triple $(E,\phi',\lambda)$ forms a rank $r$ Drinfeld $A$-module with level $tn$-structure over $\Spec(RS_{W,0})$.

\begin{Thm}\label{moduli}
The affine scheme $M^r_{tn,B}:=\Spec(RS_{W,0})$, together with the universal family $(E,\phi',\lambda)$, represents the functor from $B$-Schemes to Sets which sends a scheme $S$ over $\Spec(B)$ to the set of isomorphism classes of triples $(L,\psi,\mu)_S$, where $(L,\psi)$ is a rank $r$ Drinfeld $A$-module over $S$, and $\mu : W \to L(S)$ is a level-$tn$ structure.  
\end{Thm}

The special case $M^r_{t,B}\cong\Spec(RS_{V,0})$ is essentially due to Pink \cite[\S 7]{Pink} 
and inspired Theorem~\ref{moduli}.

\medskip

\begin{Proof}
Let $S$ be a scheme over $\Spec(B)$ and $(L,\psi,\mu)_S$ a triple as above. 
We must associate to the isomorphism class of $(L,\psi,\mu)_S$ an $S$-valued point $\eta$ on $\Spec(RS_{W,0})$ such that the pullback of the universal family $(E,\phi',\lambda)$ to $\eta$ is isomorphic to $(L,\psi,\mu)_S$.


First notice that the line bundle $L/S$ must be trivial, since for any $v\in\Vo$, $\mu(v)\in L(S)$ is a nowhere zero section, as $t$ is prime to the characteristic of $\psi$. Now cover $S$ with open affines $\Spec(R)$; it suffices to prove that the isomorphism class of each pullback $(L,\psi,\mu)_{\Spec(R)}$ corresponds to a $\Spec(R)$-valued point on $\Spec(RS_{W,0})$. Thus we assume that $S=\Spec(R)$ is affine, where $R$ is a $B$-algebra, and that $L=\BG_{\mathrm{a},R}$ is the additive group scheme over $\Spec(R)$.

Next, we may replace $\psi$ by an isomorphic Drinfeld module which is standard, i.e. for which
\[
\psi_t(X) = tX+a_1X^q+\cdots+a_rX^{q^r},
\]
where $a_1,\ldots,a_{r-1}\in R$ and $a_r\in R^*$, see \cite[\S2.2.3, p21]{Lehmkuhl}.

The level structure $\mu$ is a morphism $\mu : W \longto R$ such that $\mu(\Wo)\subset R^*$ 
and
\begin{eqnarray*}
\psi_t(X) & = & tX\prod_{v\in\Vo}\left(1-\frac{X}{\mu(v)}\right), \\
\psi_{tn}(X) & = & tnX\prod_{w\in\Wo}\left(1-\frac{X}{\mu(w)}\right).
\end{eqnarray*}

Recall that we have fixed $v_1\in\Vo$. Consider the Drinfeld module $\psi' := \mu(v_1)^{-1}\psi\mu(v_1)$; it is isomorphic to $\psi$ over $R$, and 
\begin{eqnarray*}
\psi'_t(X) & = & tX\prod_{v\in\Vo}\left(1-\frac{\mu(v_1)}{\mu(v)}X\right), \\ 
\psi'_{tn}(X) & = & tnX\prod_{w\in\Wo}\left(1-\frac{\mu(v_1)}{\mu(w)}X\right). 
\end{eqnarray*}

We now consider the $B$-algebra homomorphism
\[
\theta : RS_{V,0} \longto R
\]
which is determined by 
\[
\frac{v}{v'} \longmapsto \frac{\mu(v)}{\mu(v')}, \qquad \text{for $v,v'\in\Vo$.}
\]
This exists because $M^r_{t,B}\cong\Spec RS_{V,0}$, by \cite[\S 7]{Pink}, but one can also see this directly: all relations satisfied by the $v/v'$ in $RS_{V,0}$ are also satisfied by the $\mu(v)/\mu(v')$ in $R$.
 
We extend $\theta$ to the $B$-algebra homomorphism
\[
\theta : RS_{V,0}[w'_1,w'_2,\ldots,w'_r] \longto R
\]
determined by
\begin{eqnarray*}
w'_i & \longmapsto & \frac{\mu(w_i)}{\mu(v_1)}. 
\end{eqnarray*}
It is clear that $I_{W,0}\subset \ker\theta$, so $\theta$ extends to a $B$-algebra homomorphism $\theta : RS_{W,0} \longto R$. 
Furthermore, $\theta$ defines a $\Spec(R)$-valued point $\eta$ on $\Spec(RS_{W,0})$, and the triple $(L,\psi,\mu)$ is isomorphic to the pullback $(L,\psi',\mu(v_1)^{-1}\mu)$ of the universal family $(E,\phi',\lambda)$ to $\eta$, as required.

\medskip

Conversely, suppose we are given an $S$-valued point $\eta$ on $\Spec(RS_{W,0})$,
then the pullback of the universal family $(E,\phi',\lambda)$ along $\eta : S \to \Spec(RS_{W,0})$ defines a triple with the desired properties.
\end{Proof}

Next, we collect here the following fundamental results on Drinfeld moduli schemes.

\begin{Thm}\label{Drinfeld}
Let $n\in A=\BF_q[t]$ be monic and recall that $B=A[\frac{1}{tn}]$.
\begin{enumerate}
  \item[(i)] The scheme $M_{tn,B}^r$ is smooth of relative dimension $r-1$ over $\Spec B$.
	\item[(ii)] The group $\GL_r(A/tnA)$ acts on the level structure of the universal Drinfeld module $\phi'$, and this induces an action of $\GL_r(A/tnA)$ on $M_{tn,B}^r$.
  \item[(iii)] The canonical morphism $M_{tn,B}^r\to M_{t,B}^r$ is \'etale with Galois group 
	\[
	G_r(n) := \ker\big(\GL_r(A/tnA)\longto\GL_r(A/tA)\big).
	\]
	\item[(iv)] There is a morphism, defined over $\Spec B$, 
	\[
	w_{tn} : M_{tn,B}^r \longto M_{tn,B}^1,
	\]
	which is compatible with the action of $\GL_r(A/tnA)$, in the sense that, for every $\sigma\in\GL_r(A/tnA)$, we have $w_{tn}\circ\sigma = \det(\sigma)\circ w_{tn}$. 
	\item[(v)] The scheme $M_{tn,B}^r$ is integral, and the rings $RS_{W,0}$ and $RS_W$ are integral.
\end{enumerate}
\end{Thm}

\begin{Proof}
The first three statements are essentially due to Drinfel'd \cite{Drinfeld}, who proved this more generally over $\Spec A$ but for level structures divisible by two distinct primes. In our situation, the level $tn\neq 1$ is invertible in $B$. Thus, as in the proof of Theorem \ref{moduli}, if $(L,\psi,\mu)_S$ is a Drinfeld module with level-$tn$ structure over a $B$-scheme $S$ then any $v\in\Vo$ gives a nowhere vanishing section $\mu(v)\in L(S)$, which trivializes $L/S$.  For $S$ over $\Spec A$ such a trivialization is only achieved if the level structure is divisible by two distinct primes, see \cite[Prop. 2.5.1 and Theorem 3.4.1]{Lehmkuhl} for details.  Thus in our case, Drinfeld's proofs give (i), (ii) and (iii) above. See also \cite{vdHeidenMC}, as well as \cite{Hubschmid} for a very clear exposition of the situation over the quotient field of $A$. 

Alternatively, the interested reader is challenged to deduce (i)--(iii) directly from Theorem~\ref{moduli}, for example the fact that $\Spec(RS_{W,0})\to\Spec(RS_{V,0})$ is \'etale follows exactly as in Proposition \ref{Prop:rings}. 

Statement (iv) is essentially due to Anderson \cite{Anderson}, see also \cite{vdHeidenMC} for details.

To prove (v), note that $RS_{W,0}$ is flat over $B$, by (i), so $RS_{W,0}$ injects into $RS_{W,0}\otimes_B F$, which is integral by \cite[Cor. 3.4.5]{Hubschmid}. 
%
%
%
%

Lastly, $RS_W = RS_{W,0}[v_1]$, and $v_1$ is transcendental over $RS_{W,0}$, so $RS_W$ is also integral. 
\end{Proof}

\section{Sufficiently generic Drinfeld modules}\label{Sec:generic}

%

Now let $\psi$ be a rank $r$ Drinfeld $A$-module over an integral $B$-algebra $R$ defined by
\[
\psi_t(X) = tX + a_1X^q + \cdots + a_rX^{q^r},
\]
where $a_1,\ldots,a_{r-1}\in R$ and $a_r\in R^*$. We define the invariants
\[
J_i := \frac{a_i^{(q^r-1)/d_i}}{a_r^{(q^i-1)/d_i}}, \qquad i=1,\ldots,r-1,
\]
where $d_i := \gcd(q^i-1,q^r-1)$. (Actually, we could choose $d_i$ to be any common divisor of $q^i-1$ and $q^r-1$.) These are isomorphism invariants, although for $r\geq 3$ they do not determine the isomorphism class of $\psi$ completely, see \cite{Potemine}.  

\begin{Def}
A Drinfeld module $\psi$ of rank $r\geq 1$ is {\em sufficiently generic} if $r=1$, or if $r\geq 2$ and the invariants $J_1,\ldots,J_{r-1}$ are algebraically independent over $\BF_q(t)$.
\end{Def}

This condition is equivalent to the ring of isomorphism invariants (see \cite{Potemine}) of $\psi$ having transcendence degree $r$ over $\BF_q$.

Consider the subfield $K:=\BF_q(t,a_1,\ldots,a_r)$ of the quotient field of $R$, and denote by $K_{tn}$ the splitting field of $\psi_{tn}(X)$ over $K$. We denote by 
\[
RS_{tn,0}:=B\big[\frac{w}{w'} \;|\; w,w'\in\psi[tn], w'\neq 0\big]
\] 
the $B$-subalgebra of $K_{tn}$ generated by the quotients $\frac{w}{w'}$ with $w,w'\in\psi[tn]$, $w'\neq 0$. 

Our second main result is the following.

\begin{Thm}\label{generic}
If $\psi$ is sufficiently generic, then $RS_{tn,0}\cong RS_{W,0}$. In particular, 
\[
M^r_{tn,B} \cong \Spec(RS_{tn,0}).
\]
\end{Thm}

\begin{Proof}
When $r=\rank(\psi)=1$ we can show directly that $RS_{tn,0} \cong RS_{W,0}$. Let $u_1$ be a generator of $\psi[t]\cong A/tA$, and set $\psi':=u_1^{-1}\psi u_1$. Then 
\[
\psi'_t(X) = u_1^{-1}\psi_1(u_1X) = tX\prod_{\varepsilon\in\BF_q^*}\left(1-\frac{u_1X}{\varepsilon u_1}\right) = \phi'_t(X).
\]
Now 
\begin{eqnarray*}
RS_{tn,0} & =  & B\big[\frac{w}{u_1}, \frac{u_1}{w} \;|\; 0\neq w\in\psi[tn]\big] = B\big[w', \frac{1}{w'} \;|\; 0\neq w'\in\psi'[tn]\big] \\
& \cong & B\big[w', \frac{1}{w'} \;|\; 0\neq w'\in\phi'[tn]\big]. 
\end{eqnarray*} 
But the last expression is equal to $B\big[\frac{w}{w'} \;|\; w,w'\in\phi[tn]\smallsetminus\{0\}\big]\cong RS_{W,0}$.

Now suppose that $r\geq 2$.
Choose a level-$tn$ structure $\mu : W \to \psi[tn]\subset K_{tn}$. Then $\mu(\Wo)\subset K_{tn}^*$ and 
similarly to part 1 of the proof of Theorem \ref{moduli}, we construct a $B$-algebra homomorphism
\[
\theta : RS_W \longto K_{tn}; \qquad \theta(w_i)=\mu(w_i), \; \theta(z)=\prod_{w\in\Wo}\mu(w)^{-1}\in K_{tn}^*, \qquad i=1,2,\ldots,r.
\]

We must show that $\ker \theta \cap RS_{W,0} = \{0\}$, so suppose that $0\neq f \in \ker\theta\cap RS_{W,0}$. By Theorem~\ref{Drinfeld}.(iii), $\prod_{\sigma\in G_r(n)}\sigma(f)\in RS_{V,0}$. Multiplying this by a suitable unit 
$u\in RS_V^*$,
we obtain a homogeneous element
\[
\tilde{f} = u \prod_{\sigma\in G_r(n)}\sigma(f) \in \ker\theta\cap R_V.
\]
Now, by \cite[Theorem 3.1]{PinkSchieder}, $\GL_r(\BF_q)$ acts on $R_V$ and the ring of invariants is $R_V^{\GL_r(\BF_q)} = B[g_1,\ldots,g_r]$, where
\[
\phi_t(X) = tX\prod_{v\in\Vo}\left(1-\frac{1}{v}X\right) = tX + g_1X^q + \cdots + g_rX^{q^r}.
\]
Thus we obtain
\[
\bar{f} := \prod_{\tau\in\GL_r(\BF_q)}\tau(\tilde{f})\in \ker\theta\cap B[g_1,\ldots,g_r]
\]
which is homogeneous of some degree $d$ with respect to the grading $\deg(g_i)=q^i-1$ for $i=1,\ldots,r$. Notice that $\bar{f}\neq 0$ since $RS_{W}$ is integral, by Theorem \ref{Drinfeld}.(v).

Since $a_i=\theta(g_i)$ for $i=1,\ldots,r$, we see that
\[
\bar{f}(a_1,\ldots,a_r)=0.
\]

Now, let $\delta\in \Ksep$ be such that $\delta^{q^r-1}=a_r$, and set 
\[
u_i := \delta^{1-q^i}a_i,\quad\text{so}\quad J_i=u_i^{(q^r-1)/d_i},   \qquad i=1,\ldots,r-1.
\]
Then $0=\delta^d\bar{f}(u_1,u_2,\ldots,u_{r-1},1)$. It follows that $\BF_q(t,u_1,\ldots,u_{r-1})$ has transcendence degree at most $r-2$ over $\BF_q(t)$. Since $\BF_q(t,J_1,\ldots,J_{r-1})\subset \BF_q(t,u_1,\ldots,u_{r-1})$, this contradicts the algebraic independence of $J_1,\ldots,J_{r-1}$ over $\BF_q(t)$.
\end{Proof}

%

\section{Abyhankar's Generalized Iteration Conjecture}\label{Sec:Abhyankar}

We now come to the heart of this article.

Let $k$ be a field containing $\BF_q$ such that $t$ is transcendental over $k$, and set $F:=k(t)$ and $K:=F(a_1,a_2,\ldots,a_{r-1})$, where $a_1,a_2,\ldots,a_{r-1}$ are algebraically independent over $F$. Consider the rank $r$ Drinfeld module $\psi$ over $K$ defined by 
\[
\psi_t(X) = tX + a_1X^q + \cdots + a_{r-1}X^{q^{r-1}} + X^{q^r}.
\]
Notice that here the highest coefficient is $a_r=1$.

Let $n\in A$ be any monic polynomial and denote by $K_n$ the splitting field of $\psi_n(X)$ over $K$. Our third main result is the following,
which was conjectured by S.~S.~Abhyankar in \cite[\S19]{Ab01}:

\begin{Thm}[Generalized Iteration Conjecture]\label{Abhyankar}
$\Gal(K_{n}/K)\cong \GL_r(A/nA)$.
\end{Thm}

Equivalently, the Galois representation attached to the $n$-torsion of $\psi$ over $K$ is surjective, for every non-zero $n\in A$.

Since $[K_n:K]$ does not depend on $k$, we obtain
\begin{Cor}
$K_n/K$ is a purely geometric extension.
\end{Cor}

A number of special cases of Theorem \ref{Abhyankar} are known, see for example \cite{AS1,Thiery} and further references in \cite[\S19]{Ab01}. A related result is in \cite{Joshi}. In particular, the case $r=1$ follows from the work of Carlitz \cite{Carlitz}, while the case where $n=t$ dates back to E.~H.~Moore \cite{Moore}:

\begin{Thm}[Moore]\label{ThmMoore}
$\Gal(K_t/K)\cong\GL_r(\BF_q)$.
\end{Thm}

\begin{Proof} 
See \cite[\S3]{AS1} for a particularly simple proof.
\end{Proof}

It is clear that $\Gal(K_n/K)$ is isomorphic to a subgroup of $\GL_r(A/nA)$, and this subgroup cannot be enlarged by enlarging $k$, so we may  
\[
\text{assume that}\quad\BF_{q^r} \subset k.
\]

Denote by $K_{tn}$ the splitting field of $\psi_{tn}(X)$ over $K$. It will be sufficient to prove the following result.

\begin{Prop}\label{relative}
Let $n\in A$ be monic. Then 
\[
\Gal(K_{tn}/K_t)\cong G_r(n) = \ker\big(\GL_r(A/tnA)\longto\GL_r(A/tA)\big).
\]
\end{Prop}

\begin{Proofof}{Theorem \ref{Abhyankar}}

\begin{minipage}[u]{5cm}
\[
\xymatrix{
& K_{tn}\ar@{-}[dl]_{G_r(n)}\ar@{-}[dr] & \\
K_t\ar@{-}[dr]_{\GL_r(A/tA)} & & K_n\ar@{-}[dl] \\
& K &
}
\]
\end{minipage}
\hspace{0.5cm}
\begin{minipage}[u]{9cm}
{By Proposition \ref{relative} $\Gal(K_{tn}/K_t)\cong G_r(n)$. Consider the field extensions in the diagram. We have $\Gal(K_t/K)\cong\GL_r(A/tA)$ by Theorem \ref{ThmMoore}, thus $[K_{tn}:K] = \#G_r(n)\cdot\#\GL_r(A/tA) = \#\GL_r(A/tnA)$ and it follows that $\Gal(K_{tn}/K)\cong\GL_r(A/tnA)$. Now from the action of $\Gal(K_{tn}/K)$ on $\psi[n]$ we see that $\Gal(K_{tn}/K_n)\cong \ker\big(\GL_r(A/tnA)\to\GL_r(A/nA)\big)$ and thus $\Gal(K_n/K)\cong\GL_r(A/nA)$.}

\end{minipage}

\end{Proofof}

It remains to prove Proposition \ref{relative}. This will be the goal of the rest of this paper.

\section{Moore and Carlitz}\label{Sec:MooreCarlitz}

We associate to $\psi$ its {\em determinant Drinfeld module} $\rho$, which is the rank 1 Drinfeld module defined over $F$ by
\[
\rho_t(X) = tX -(-1)^ra_rX^q = tX -(-1)^rX^q.
\]
When $r$ is even, then $\rho$ is the original Carlitz module (as studied by Carlitz in the 1930s, \cite{Carlitz}), whereas, when $r$ is odd then $\rho$ is the ``modern'' Carlitz module, as defined in modern texts such as \cite[Chapter 3]{GossBS} and \cite[Chapter 12]{Rosen}. 
 
We denote by $F_t$ and $F_{tn}$ the splitting fields of $\rho_t(X)$ and $\rho_{tn}(X)$ over $F$, respectively. 

\begin{Prop}[Carlitz]\label{Carlitz}
We have $\Gal(F_{tn}/F)\cong (A/tnA)^*$ and $\Gal(F_t/F)\cong \BF_q^*$.
\end{Prop}

\begin{Proof}
Let $F'=\BF_q(t)$ and denote by $F'_{tn}$ the splitting field of $\rho_{tn}(X)$ over $F'$. Denote by $\BF$ the algebraic closure of $\BF_q$ in $k$.

L.~Carlitz proved in 1938 that $\Gal(F'_{tn}/F')\cong (A/tnA)^*$ (\cite{Carlitz}, see also \cite[Theorem 12.8]{Rosen}).
Furthermore, the extension $F'_{tn}/F'$ is purely geometric, by \cite[Corollary to Theorem 12.14]{Rosen}, so also
$\Gal(\BF F'_{tn}/\BF F')\cong (A/tnA)^*$.
%

Lastly, since $t$ is transcendental over $k$, we have $F\cap (\BF F'_{tn}) = \BF F'$, and so 
\[
\Gal(F_{tn}/F) = \Gal(F \BF F'_{tn} / F \BF F') \cong \Gal(\BF F'_{tn}/ \BF F') \cong (A/tnA)^*,
\]
and $\Gal(F_t/F)\cong \BF_q^*$ follows by setting $n=1$.
%
%
%
\end{Proof}

The determinant Drinfeld module $\rho$ plays the same role for $\psi$ that the multiplicative group $\BG_{\mathrm m}$ plays for elliptic curves, and the analogue of the Weil Pairing, developed in \cite{Anderson,vdHeiden} in general, has a particularly simple description in the case of $t$-torsion using the Moore determinant. Recall (see \cite[\S1.3]{GossBS}) that the Moore determinant of a tuple $(x_1,x_2,\ldots,x_r)$ of elements in a field containing $\BF_q$ is defined by
\[
M(x_1,x_2,\ldots,x_r) := \left|\begin{array}{llll}x_1 & x_2 & \cdots & x_r \\ x_1^q & x_2^q & \cdots & x_r^q \\ \vdots & \vdots & & \vdots \\ x_1^{q^{r-1}} & x_2^{q^{r-1}} & \cdots & x_r^{q^{r-1}} \end{array}\right|
\]
and has the property that $M(x_1,x_2,\ldots,x_r)\neq 0$ if and only if $x_1,x_2,\ldots,x_r$ are linearly independent over $\BF_q$.

Choose a basis $v_1,v_2, \ldots, v_r$ of the vector space $\psi[t]\cong \BF_q^r$, then we have
\[
\psi_t(X) = M(v_1,v_2,\ldots,v_r,X)/M(v_1,v_2,\ldots,v_r),
\]
since both sides equal the unique monic polynomial with set of roots $\BF_qv_1+\BF_qv_2 + \cdots +\BF_qv_r$. Comparing $X$-coefficients gives $t = (-1)^{r}M(v_1,v_2,\ldots,v_r)^{q-1}$, so we see that $M(v_1,v_2,\ldots,v_r)\in\rho[t]$. Thus the Moore determinant defines a map (the analogue of the Weil pairing for $t$-torsion):

\[
M : (\psi[t])^r \longto \rho[t]; \qquad (x_1,x_2,\ldots,x_r) \longmapsto M(x_1,x_2,\ldots,x_r).
\]

The following result is easily verified directly.

\begin{Prop}\label{WeilPairing1}
The map $M$ above is $\BF_q$-multilinear, alternating and surjective. It follows that $F_t\subset K_t$. \qed
\end{Prop}

Via the choice of basis $v_1,v_2,\ldots,v_r$ for $\psi[t]$ we identify $\Gal(K_t/K)$ with $\GL_r(\BF_q)$, see Theorem~\ref{ThmMoore}. Since $K/F$ is purely transcendental, we also have 
\[
\Gal(KF_t/K)\cong\Gal(F_t/F)\cong\BF_q^*=\det(\GL_r(\BF_q))
\]
and
\[
\Gal(KF_{tn}/KF_t)\cong\Gal(F_{tn}/F_t)\cong G_1(n) = \ker\big((A/tnA)^* \longto (A/tA)^*\big).
\]

A direct computation shows the following.

\begin{Prop}\label{WeilPairing2}
Let $\sigma\in\Gal(K_t/K)=\GL_r(\BF_q)$ and $(x_1,x_2,\ldots,x_r)\in(\psi[t])^r$. Then 
\[
M\big(\sigma(x_1),\sigma(x_2),\ldots,\sigma(x_r)\big) = \det(\sigma)(M(x_1,x_2,\ldots,x_r)).
\]
In particular, $KF_t$ is the fixed field of $\SL_r(\BF_q)$ in $K_t$. \qed
\end{Prop}

We summarise our progress thus far in the following diagram of field extensions and Galois groups.

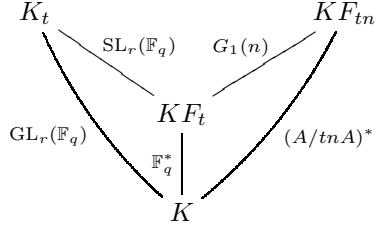
\[
\xymatrix@C+10pt{
K_t\ar@{-}@/_/[ddr]_{\GL_r(\BF_q)}\ar@{-}[dr]^{\SL_r(\BF_q)} & & KF_{tn}\ar@{-}[dl]_{G_1(n)}\ar@{-}@/^/[ddl]^{(A/tnA)^*} \\
& KF_t\ar@{-}[d]_{\BF_q^*} & \\
& K &
}
\]

\section{Function fields of Drinfeld modular varieties}\label{Sec:FF}

We define the following fields:

\begin{eqnarray*}
K_{tn,0} & := & F\left(\frac{w}{w'} \;|\; w,w'\in\psi[tn],\; w'\neq 0\right) \subset K_{tn},\\
K_{t,0} & := & F\left(\frac{w}{w'} \;|\; w,w'\in\psi[t],\; w'\neq 0\right) \subset K_{t}, \qquad\text{and} \\
F_{tn,0} & := & F\left(\frac{w}{w'} \;|\; w,w'\in\rho[tn],\; w'\neq 0\right) \subset F_{tn}.
\end{eqnarray*}

Notice that the leading coefficient of 
\[
\psi_t(X) = tX\prod_{v\in\psi[t]\smallsetminus\{0\}}\left(1-\frac{X}{v}\right)
\]
is 
\[
1 = t\prod_{v\in\phi[t]\smallsetminus\{0\}}\frac{1}{v}.
\]
Thus 
\[
v_1^{q^r-1} = t\prod_{v\in\phi[t]\smallsetminus\{0\}}\frac{v_1}{v} \in K_{t,0},
\]
and since $K_t = K_{t,0}(v_1)$ and we have assumed that $\BF_{q^r}\subset k\subset K_{t,0}$, we obtain

\begin{Prop}
The extension $K_t/K_{t,0}$ is Galois with $C:=\Gal(K_t/K_{t,0})$ cyclic of order dividing $q^r-1$. \qed
\end{Prop}

\begin{Rem}
With a little more effort one can show that in fact $C$ has order equal to $q^r-1$, but we will not need this here.
\end{Rem}

We have 
\[
K_{t,0}=F\big(\frac{v_2}{v_1},\frac{v_3}{v_1},\ldots,\frac{v_r}{v_1}\big), 
\]
for any basis $v_1,v_2,\ldots,v_r$ of $\psi[t]$, 
and since $K_{t,0}$ has transcendence degree $r$ over $\BF_q$ (because $K_t$ does) it follows that
$K_{t,0}/F$ is a purely transcendental extension of transcendence degree $r-1$.

Furthermore, since $K_t$ contains a generator of $\rho[t]\subset\rho[tn]$, we see that $K_tF_{tn,0}=K_tF_{tn}$.

\begin{Prop}\label{Prop:FF}
We have
\begin{enumerate}
	\item[(i)] $\Gal(K_{tn,0}/K_{t,0}) \cong G_r(n) = \ker\big(\GL_r(A/tnA)\longto\GL_r(A/tA)\big)$.
	\item[(ii)] The subfield of $K_{tn,0}$ fixed by
	\[
	S_r(n) := \ker\big( \SL_r(A/tnA) \longto \SL_r(A/tA) \big)
	\]
	is $K_{t,0}F_{tn,0}$.
\end{enumerate}
\end{Prop}

\begin{Proof}
We first consider the special case where $k=\BF_q$ and $F=\BF_q(t)$. 

The base extension $M_{tn,F}^r = M^r_{tn,B}\times_{\Spec B}\Spec F$ is integral, by Theorem \ref{Drinfeld}.(v), and, since $\psi$ is sufficiently generic, its function field over $F$ is $K_{tn,0}$, by Theorem \ref{generic}. Similarly, the function fields of $M^r_{t,F}$, $M^1_{tn,F}$ and $M^1_{t,F}$ over $F$ are $K_{t,0}$, $F_{tn,0}$ and $F_{t,0}=F$, respectively.
Now (i) follows from Theorem \ref{Drinfeld}.(iii).

To prove (ii), the fixed field contains $K_{t,0}F_{tn,0}$, by Theorem \ref{Drinfeld}.(iv), while $\Gal(K_{t,0}F_{tn,0}/K_{t,0})\cong \Gal(F_{tn,0}/F)$, since $K_{t,0}$ is purely transcendental over $F_{t,0}=F$. Now $\Gal(F_{tn,0}/F)\cong G_1(n)$ (by Theorem \ref{Drinfeld}.(iii)), which is isomorphic to the quotient $G_r(n)/S_r(n)$. The result follows in this case.

To extend our result to the case for general $k$, recall that $t,a_1,\ldots,a_{r-1}$ are algebraically independent over $k$, so it suffices to show that the relevant field extensions are purely geometric, i.e. that $\BF_q$ is algebraically closed in the function field of $M^r_{tn,\BF_q(t)}$ over $\BF_q(t)$. We achieve this by constructing a field $L$, in which $\BF_q$ is algebraically closed, and a rank $r$ Drinfeld $\BF_q[t]$-module $\rho'$ over $L$ with $\rho'[tn]\subset L$.

Let $A'=\BF_q[\sqrt[r]{t}]$ and $K'=\BF_q(\sqrt[r]{t})$. Consider the Carlitz $A'$-module $\rho'$ defined over $K'$ by 
\[
\rho'_{\sqrt[r]{t}}(X) = \sqrt[r]{t}X+X^q.
\]
As before, $L:=K'(\rho'[tn])$ is purely geometric over $K'$. On the other hand, $\rho'$ is also a rank $r$ Drinfeld $\BF_q[t]$-module (with complex multiplication by $A'$), so it, together with a level-$tn$ structure over $L$, defines an $\BF_q(t)$-algebra homomorphism $RS_{W,0}\otimes_B \BF_q(t) \to L$. It follows that $\BF_q$ is algebraically closed in the function field of $M^r_{tn,\BF_q(t)}$ over $\BF_q(t)$.
\end{Proof}

We summarise our progress in the following diagram.

\[
\xymatrix{
K_{tn,0}\ar@{-}[ddr]_{G_r(n)}\ar@{-}[dr]^{S_r(n)} & & K_tF_{tn,0} = K_tF_{tn}\ar@{-}[dl]\ar@{-}[d] \\
& K_{t,0}F_{tn,0}\ar@{-}[d]^{G_1(n)} & K_t\ar@{-}[dl]^{C} \\
& K_{t,0} &
}
\]

Since the order of $C$ is prime to $p$, we see that 

\begin{Prop}\label{p-order1}
We have $v_p \big([K_tF_{tn}:K_t]\big) = v_p \big(|G_1(n)|\big),$ where $v_p$ denotes the $p$-adic valuation.\qed
\end{Prop}

\section{Some Group Theory}\label{Sec:groups}

Before we continue, we need to recall some results from group theory.

\begin{Lem}\label{Groups1}
Every proper Abelian quotient of $\SL_r(\BF_q)$ has order $p$.
\end{Lem}

\begin{Proof}
If we use $\der$ to denote the derived (commutator) subgroup, then by \cite[chap. XIII Theorems 8.3 and 9.2]{LangAlg} we have
\[
\SL_r(\BF_q)^\der = \SL_r(\BF_q),
\]
with two exceptions. These are:
\begin{itemize}
	\item If $r=2$ and $q=2$, then $\SL_2(\BF_2)^\der \cong A_3$, which has index 2 in $\SL_2(\BF_2)\cong S_3$, and
	\item If $r=2$ and $q=3$, then $\SL_2(\BF_3)^\der \cong Q$, the 8-element quaternion group, which has index 3 in $\SL_2(\BF_3)$. 
\end{itemize}
The result follows.
\end{Proof}

\begin{Prop}\label{Groups2}
Every proper Abelian quotient of $S_r(n)$ is a $p$-group.
\end{Prop}

\begin{Proof}
Let the prime factorisation of $n$ in $A$ be given by
\[
n = \prod_P P^{a_P}.
\]
Then
\begin{eqnarray*}
S_r(n) & = & \ker\big( \SL_r(A/tnA) \longto \SL_r(A/tA) \big) \\
& \cong & \ker\big( \SL_r(A/t^{a_t+1}A) \longto \SL_r(A/tA) \big) \times \prod_{P|n, P\neq t} \SL_r(A/P^{a_P}A).
\end{eqnarray*}
For every prime polynomial $P\in A$, the group $\ker\big( \pi: \SL_r(A/P^{a}A) \longto \SL_r(A/PA) \big)$ is a $p$-group (of order $q^{\deg(P)(a-1)(r^2-1)}$).

It remains to show that any Abelian quotient of the form $\SL_r(A/P^aA)/N$ is a $p$-group. Write $\#\SL_r(A/P^aA) = p^bm$, where $p\nmid m$. Since $\pi(N) < \SL_r(A/PA)^\der$, 
$\SL_r(A/PA)/\pi(N)$ is a $p$-group by Lemma \ref{Groups1}, thus $m|\#\pi(N)$ and so also $m|\# N$, since $\ker(\pi)$ is a $p$-group. The result follows.
%
%
%
%
%
%
\end{Proof}

\section{Completing the proof}\label{Sec:proof}

We now have all the ingredients we need. Our next step is

\begin{Prop}\label{Intersection1}
$K_t\cap KF_{tn} = KF_t$. In particular, 
\[
 \Gal(K_tF_{tn}/K_t)\cong G_1(n) \qquad\text{and}\qquad  \Gal(K_tF_{tn}/KF_{tn})\cong \SL_r(\BF_q).
\]
\end{Prop}

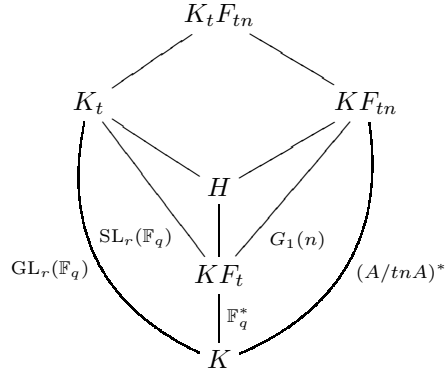
\begin{Proof}
\[
\xymatrix@R-5pt{
& K_tF_{tn}\ar@{-}[dl]\ar@{-}[dr] & \\
K_t\ar@{-}[dr]\ar@{-}[ddr]_>>>>>>{\SL_r(\BF_q)}\ar@{-}@/_2pc/[dddr]_{\GL_r(\BF_q)} & & KF_{tn}\ar@{-}[dl]\ar@{-}[ddl]^>>>>>>{G_1(n)}\ar@{-}@/^2pc/[dddl]^{(A/tnA)^*} \\
& H\ar@{-}[d] & \\
& KF_t\ar@{-}[d]^{\BF_q^*} & \\
& K & 
}
\]
Let $H:=K_t\cap KF_{tn}$. 
First notice that $KF_{tn}/H$ is an Abelian extension corresponding to a subgroup of $\Gal(KF_{tn}/KF_t)\cong G_1(n)$. By Proposition \ref{p-order1}, we see that $v_p\big(|G_1(n)|\big) = v_p \big([K_tF_{tn}:K_t]\big) = v_p \big([KF_{tn} : H]\big)$, and so $p\nmid [H:KF_t]$. Now $\Gal(H/KF_t)$ is Abelian of order prime to $p$; it is also a quotient of $\Gal(K_t/KF_t)\cong\SL_r(\BF_q)$, hence by Lemma \ref{Groups1} it must be trivial. The result follows.
\end{Proof}

Since we now know that 
\[
\Gal(K_tF_{tn,0}/K_t)=\Gal(K_tF_{tn}/K_t)\cong G_1(n),
\]
we see that $\Gal(K_tF_{tn,0}/K_{t,0}F_{tn,0})$ is Abelian of order prime to  $p$. Proceeding as in the proof of Proposition \ref{Intersection1}, we see that $\Gal(K_{tn,0}\cap K_tF_{tn,0} / K_{t,0}F_{tn,0})$ is an Abelian quotient of $S_r(n)$ of order prime to $p$, and hence, by Proposition \ref{Groups2}, trivial. It follows that 
\[
\Gal(K_{tn,0}K_tF_{tn,0} / K_tF_{tn,0}) \cong \Gal(K_{tn,0}/K_{t,0}F_{tn,0})\cong S_r(n),
\]
and so 
\[
\Gal(K_{tn,0}K_tF_{tn,0} / K_t) \cong G_r(n).
\]
Lastly, $K_{tn} = K_{tn,0}(v_1) \subset K_{tn,0}K_t$, so in fact $K_{tn,0}K_tF_{tn,0} = K_{tn}$ and the proof is complete. \qed
%

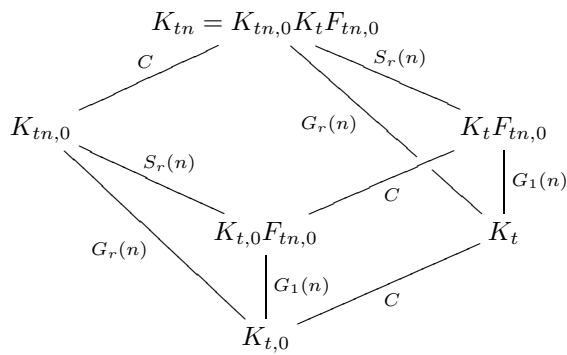
\[
\xymatrix{
& K_{tn}=K_{tn,0}K_tF_{tn,0}\ar@{-}[dl]_C\ar@{-}[dr]^{S_r(n)}\ar@{-}[ddr]_(.4){G_r(n)}|!{[dr];[dd]}\hole & \\
K_{tn,0}\ar@{-}[dr]^{S_r(n)}\ar@{-}[ddr]_{G_r(n)} & & K_tF_{tn,0}\ar@{-}[dl]^C\ar@{-}[d]^{G_1(n)} \\
& K_{t,0}F_{tn,0}\ar@{-}[d]^{G_1(n)} & K_t\ar@{-}[dl]^C \\
& K_{t,0} & 
}
\]

\begin{Rem}
Given our explicit description for the various fields concerned, it is tempting to search for a direct proof that $[K_{tn}:K_{tn,0}] =[K_{t}:K_{t,0}]$, which would allow us to cut short much of the above argument and simplify the proof of Proposition \ref{relative}. Alas, the author was not successful with this.
\end{Rem}

%



\begin{thebibliography}{99}



\bibitem{Ab01} S.~S.~Abhyankar, Resolution of singularities and modular Galois theory. {\em Bull. Amer. Math. Soc. (N.S.)} {\bf 38} (2001), no. 2, 131--169. 




\bibitem{AS1} S.~S.~Abhyankar and G.~S.~Sundaram, Galois theory of Moore-Carlitz-Drinfeld modules. {\em C. R. Acad. Sci. Paris S\'er. I Math.} {\bf 325} (1997), no. 4, 349--353. 


\bibitem{Anderson} G.~Anderson, $t$-motives. {\em Duke Math. J.} {\bf 53} (1986), no. 2, 457--502.



\bibitem{Carlitz} L.~Carlitz, A class of polynomials. {\em Trans. Amer. Math. Soc.} {\bf 43} (1938), no. 2, 167--182.

\bibitem{Drinfeld}V. G.~Drinfel'd, Elliptic modules (Russian), {\em Math. Sbornik},
{\bf 94} (1974), 594-627. Translated in {\em Math. USSR. S.}, {\bf 23} (1974), 561--
592.





\bibitem{GossBS} D.~Goss, Basic structures in function field arithmetic,
Springer-Verlag, 1996.
%
%
%




\bibitem{Hubschmid} P.~Hubschmid, The Andr\'e-Oort conjecture for Drinfeld modular varieties, {\em Compos. Math.} {\bf 149} (2013), no. 4, 507--567.

\bibitem{Joshi} K.~Joshi, A family of \'etale coverings of the affine line, {\em J. Number Theory} {\bf 59} (1996), 414--418.

\bibitem{LangEF} S.~Lang, Elliptic Functions, 2nd edition, {\em Graduate Texts in Mathematics} {\bf 112}, Springer-Verlag, 1987.

\bibitem{LangAlg} S.~Lang, Algebra, 3rd edition, {\em Graduate Texts in Mathematics} {\bf 211}, Springer-Verlag, 2002.

\bibitem{Lehmkuhl} T. Lehmkuhl, Compactification of the Drinfeld Modular Surfaces,  {\em Mem. Amer. Math. Soc.} {\bf 197} (2009), no. 921.
 
\bibitem{Moore} E.~H.~Moore, A two-fold generalization of Fermat's theorem. {\em Bull. Amer. Math. Soc.} {\bf 2} (1896), no. 7, 189--199.



\bibitem{Pink}R.~Pink, Compactification of Drinfeld modular varieties and Drinfeld Modular Forms of Arbitrary Rank, {\em Manuscripta Math.} {\bf 140} (2013), no. 3-4, 333--361.

\bibitem{PinkSchieder}R.~Pink, S.~Schieder, Compactification of a Drinfeld Period Domain over a Finite Field, {\em J. Algebraic Geometry}  {\bf 23} (2014), no. 2, 201--243.  

\bibitem{Potemine}I.~Y.~Potemine, Minimal terminal $\BQ$-factorial models of Drinfeld
coarse moduli schemes, {\em Math. Phys. Anal. Geom.} {\bf 1} (1998), 171--191.

\bibitem{Rosen} M.~Rosen, Number Theory in Function Fields, {\em Graduate Texts in Mathematics} {\bf 210}, Springer-Verlag, 2002.

\bibitem{StacksProject} The Stacks Project Authors, {\em Stacks Project}, \url{http://stacks.math.columbia.edu}, 2015.

\bibitem{Thiery} A.~Thiery, $\BF_q$-linear Galois theory, {\em J. London Math. Soc.} (2) {\bf 53} (1996), 441--454.

\bibitem{vdHeiden} G.-J.~van der Heiden, Weil pairing for Drinfeld modules, {\em Monatsh. Math.} {\bf 143} (2004), 115--143.

\bibitem{vdHeidenMC} G.-J. van der Heiden, Drinfeld modular curves and the Weil pairing, {\em J. Algebra} {\bf 299} (2006), 374--418.


\end{thebibliography}
\end{document}